\documentclass[a4paper,11pt]{amsart}

\usepackage{amsmath,amsthm,amssymb}
\usepackage[mathscr]{eucal}
\usepackage[bookmarks=false]{hyperref}

\usepackage[all]{xy}

\numberwithin{equation}{section}

\theoremstyle{plain}
\newtheorem{theorem}{Theorem}[section]

\newtheorem{lemma}[theorem]{Lemma}

\newtheorem{corollary}[theorem]{Corollary}
\theoremstyle{definition}
\newtheorem{definition}[theorem]{Definition}
\newtheorem{example}[theorem]{Example}
\newtheorem{examples}[theorem]{Examples}

\newtheorem{remark}[theorem]{Remark}

\newcommand{\N}{\mathbb{N}}

\newcommand{\C}{\mathbb{C}}

\newcommand{\F}{\mathbb{F}}
\newcommand{\cA}{\mathcal{A}}
\newcommand{\cB}{\mathcal{B}}

\newcommand{\cM}{\mathcal{M}}
\newcommand{\cN}{\mathcal{N}}

\newcommand{\cP}{\mathcal{P}}

\newcommand{\cU}{\mathcal{U}}
\newcommand{\cV}{\mathcal{V}}

\newcommand{\cZ}{\mathcal{Z}}

\newcommand{\wV}{\widetilde V}

\def\G{\Gamma}

\def\ms{\mathscr}

\begin{document}

\title[ II$_1$ factors with non-isomorphic ultrapowers]{II$_1$ factors with non-isomorphic ultrapowers}

\author[R. Boutonnet]{R\'{e}mi Boutonnet}
\address{Department of Mathematics, University of California San Diego, 9500 Gilman Drive, La Jolla, CA 92093, USA}
\email{rboutonnet@ucsd.edu}
\thanks{R. B. was supported in part by NSF  Grant DMS \#1161047, NSF Career Grant DMS \#1253402 and ANR Grant NEUMANN}

\author[I. Chifan]{Ionu\c{t} Chifan}
\address{Department of Mathematics, The University of Iowa, 14 MacLean Hall, Iowa City, IA  
52242, USA and IMAR, Bucharest, Romania}
\email{ionut-chifan@uiowa.edu}
\thanks{I. C. was supported in part by the NSF Grants DMS \#1301370 and DMS \#1600688.}

\author[A. Ioana]{Adrian Ioana}
\address{Department of Mathematics, University of California San Diego, 9500 Gilman Drive, La Jolla, CA 92093, USA, and IMAR, Bucharest, Romania}
\email{aioana@ucsd.edu}
\thanks{A.I. was supported in part by  NSF  Grant DMS \#1161047,  NSF Career Grant DMS \#1253402 and a Sloan Foundation Fellowship.}

\begin{abstract}
We prove that there exist uncountably many separable II$_1$ factors whose ultrapowers (with respect to arbitrary ultrafilters) are non-isomorphic. In fact, we prove that the families of non-isomorphic II$_1$ factors 
originally introduced by McDuff \cite{MD69a,MD69b} are such examples.
This entails the existence of a continuum of non-elementarily equivalent II$_1$ factors, thus settling a well-known open problem in the continuous model theory of operator algebras. \end{abstract}

\maketitle

\section{Introduction}
The ultrapower construction for II$_1$ factors,
originally introduced in \cite{Wr54, Sa62}, first came to prominence following McDuff's work \cite{MD69c}. It has since played a fundamental role in the study of von Neumann algebras. In particular, the analysis of ultrapowers of II$_1$ factors was a crucial ingredient in Connes' celebrated classification of injective factors \cite{Co75}. In the same paper, ultrapowers were used by Connes to formulate his famous (still unsolved) {\it embedding problem}. More recently, ultrapower techniques have been instrumental in the advances made in the classification of II$_1$ factors by Popa's deformation/rigidity theory (see e.g. \cite{Po04}). For more history on ultrapowers and ultraproducts of von Neumann algebras, see \cite{AH13}.

While ultrapowers of II$_1$ factors have been extremely useful in various applications, the following intrinsic problem remained open: how many ultrapowers (with respect to  a fixed ultrafilter) of separable II$_1$ factors exist, up to isomorphism?
Recently, a closely related problem has been considered in the emerging field of continuous model theory of operator algebras \cite{FHS09,FHS10,FHS11}: how many elementary equivalence classes of II$_1$ factors exist? 
This problem has received a lot of attention (see e.g.\ \cite{FHS11, GS14} and the survey \cite{Fa14}). The connection between the two problems stems from the continuous version of the Keisler-Shelah theorem \cite{Ke61, Sh71}. This asserts that two II$_1$ factors $M$ and $N$ are elementarily equivalent if and only if they have isomorphic ultrapowers, $M^{\cU}\cong N^{\cV}$, with respect to ultrafilters $\cU$ and $\cV$ on arbitrarily large sets.

At present, only three different elementary equivalence classes of II$_1$ factors appear in the literature.  More precisely, it was noticed in \cite{FGL06, FHS11} that, for separable II$_1$ factors,  property Gamma and the property of being McDuff are elementary properties (i.e. they are remembered by ultrapowers).  
Thus, the hyperfinite II$_1$ factor $R$, the free group factor $L(\mathbb F_2)$, and any non-McDuff separable II$_1$ factor that has property Gamma (see \cite{DL69}), are not elementarily equivalent.

By contrast, the existence of uncountably many non-isomorphic separable II$_1$ factors has been known for a long time \cite{MD69b, Sa69}.
This situation is partially explained by the fact that elementary equivalence of II$_1$ factors is a much coarser notion of equivalence than isomorphism.  An illuminating explanation of this fact is provided by a result in \cite{FHS11} which states that  any II$_1$ factor is elementarily equivalent to uncountably many non-isomorphic II$_1$ factors.

In this paper we solve the above problems, by proving the existence of a continuum of separable II$_1$ factors whose ultrapowers, with respect to any ultrafilters, are  non-isomorphic.

\subsection{Construction and statement of the result}
\label{sectionconstruction}

Our examples of II$_1$ factors with non-isomorphic ultrapowers come from McDuff's work \cite{MD69a,MD69b}. The construction relies on two functors $T_0$ and $T_1$, from the category of countable groups to itself, defined as follows.

Consider a countable group $\Gamma$. Let $\Gamma_i$, $i\geq 1,$ be isomorphic copies of $\Gamma$, and $\Lambda_i$, $i\geq 1,$ be isomorphic copies of $\mathbb Z$. We define $\widetilde\Gamma = \bigoplus_{i\geq 1}\Gamma_i$ and denote by $S_{\infty}$ the group of finite permutations of $\{1,2,...\}$. We consider the semidirect product $\widetilde\Gamma\rtimes S_{\infty}$ associated to  the action of $S_{\infty}$ on $\widetilde\Gamma$ which permutes the copies of $\Gamma$.  Following \cite{MD69b},
 \begin{itemize} 
 \item we define  $T_0(\Gamma)$ as the group generated by $\widetilde\Gamma$ and $\Lambda_i, i\geq 1,$ with the only relations that $\Gamma_i$ and $\Lambda_j$ commute for every  $i\geq j\geq 1$;
\item we define $T_1(\Gamma)$ as the group generated by $\widetilde\Gamma\rtimes S_{\infty}$ and $\Lambda_i, i\geq 1,$ with the only relations that $\Gamma_i$ and $\Lambda_j$ commute for every  $i\geq j\geq 1$.
\end{itemize}
The definition of $T_0$ is due to Dixmier and Lance in \cite[\S 21]{DL69}, who were inspired by a construction in \cite{MvN43}. 

The identification $\Gamma = \Gamma_1$ gives an embedding of $\Gamma$ inside $T_\alpha(\Gamma)$, for $\alpha \in\{ 0,1\}$. Moreover, every inclusion $\Sigma \subset\Sigma'$ of countable groups canonically induces an inclusion $T_{\alpha}(\Sigma)\subset T_{\alpha}(\Sigma')$. Hence, any sequence $\pmb\alpha = (\alpha_n)_{n \geq 1}$ of $0$'s and $1$'s, gives rise to a sequence of inclusions
\[\Gamma \subset T_{\alpha_1}(\Gamma) \subset T_{\alpha_1} \circ T_{\alpha_2}(\Gamma) \subset  T_{\alpha_1} \circ  T_{\alpha_2} \circ  T_{\alpha_3}(\Gamma) \subset \cdots\]

\begin{definition}
\label{defintro}
Given a sequence $\pmb\alpha$ of $0$'s and $1$'s we define
\begin{itemize}
\item $T_{\pmb\alpha}(\Gamma) := \Gamma$, if $\pmb\alpha$ is the empty sequence;
\item $T_{\pmb\alpha}(\Gamma) := T_{\alpha_1} \circ T_{\alpha_2} \circ \cdots \circ T_{\alpha_n}(\Gamma)$, if $\pmb\alpha = (\alpha_1,\alpha_2,\dots,\alpha_n)$ is a finite sequence;
\item $T_{\pmb\alpha}(\Gamma)$ is the inductive limit of the increasing sequence of groups $(T_{(\alpha_1,\dots,\alpha_n)}(\Gamma))_{n \geq 1}$, if $\pmb\alpha = (\alpha_n)_{n \geq 1}$ is an infinite sequence.
\end{itemize}
We denote by $M_{\pmb\alpha}(\Gamma) := L(T_{\pmb\alpha}(\Gamma))$ the associated von Neumann algebra.
\end{definition}

The countable family of non-isomorphic II$_1$ factors constructed in \cite{MD69a} is just $M_{\pmb\alpha_n}(\mathbb F_2)$, $n \geq 1$, where $\pmb\alpha_n$ denotes the finite $0$-valued sequence of length $n$ and $\mathbb F_2$ is the free group on two generators.
The uncountable family $M_{\pmb\alpha}(\Gamma)$, indexed over infinite sequences $\pmb\alpha$, is precisely the family of non-isomorphic II$_1$ factors constructed in \cite{MD69b}.

\begin{theorem}\label{main}
Consider a countable group $\Gamma$ and two different sequences $\pmb\alpha \neq \pmb\beta$ of $0$'s and $1$'s.
Assume that  $\Gamma=\mathbb F_2$, or that the sequences $\pmb\alpha$ and $\pmb\beta$ are infinite.

Then $M_{\pmb\alpha}(\Gamma)^\cU$ is not isomorphic to $M_{\pmb\beta}(\Gamma)^\cV$, for any ultrafilters $\cU$ and $\cV$.
\end{theorem}

We will actually deduce Theorem \ref{main} from two results.
\begin{itemize}
\item Firstly, we prove that any ultrapower of $M_{\pmb\alpha}(\F_2)$ remembers the length of $\pmb\alpha$. This is achieved by introducing a new invariant for von Neumann algebras, called {\it McDuff depth}, or {\it depth} for short, which quantifies property Gamma. We will show that if $\pmb\alpha$ is a sequence of $0$'s and $1$'s, then the depth of any ultrapower of $M_{\pmb\alpha}(\Gamma)$ is at least the length of $\pmb\alpha$, with equality for $\Gamma = \F_2$. This first half can be found in Section \ref{sectiondepth}; see Theorem \ref{countable}.
\item Secondly, we show that two infinite different sequences give different ultrapowers. For this we generalize McDuff's {\it property V} \cite{MD69b}. Using our notion of depth, we show that any  ultrapowers of $M_{\pmb\alpha}(\Gamma)$ has property $V$ at depth $k$ if and only if $\alpha_k = 1$. Hence the sequence $\pmb\alpha$ is an invariant of $M_{\pmb\alpha}(\Gamma)$ and its ultrapowers. This part is done in Section \ref{sectionpropV}; see Theorem \ref{uncountable}.
\end{itemize}

As a consequence of Theorem \ref{main} we deduce the existence a continuous family of separable non-nuclear $\cZ$-stable C$^*$-algebras with non-isomorphic ultrapowers. We thank Ilijas Farah for pointing this out to us. For any sequence $\pmb\alpha$ of $0$'s and $1$'s and any group $\Gamma$, define $A_{\pmb\alpha}(\Gamma) = C^*_r(T_{\pmb\alpha}(\Gamma)) \otimes \cZ$, where $\cZ$ is the Jiang-Su algebra.

\begin{corollary}\label{maincor}
In the setting of Theorem \ref{main}, assume moreover that $\pmb\alpha$ and $\pmb\beta$ are non-empty. 

Then the C$^*$-algebraic ultrapowers $A_{\pmb\alpha}(\Gamma)^\cU$ and $A_{\pmb\beta}(\Gamma)^\cV$ are not isomorphic, for any ultrafilters $\cU$ and $\cV$. 
\end{corollary}

Note that the proof of Corollary \ref{maincor} also implies that the reduced group C$^*$-algebras $C^*_r(T_{\pmb\alpha}(\Gamma))$ and $C^*_r(T_{\pmb\beta}(\Gamma))$ do not have isomorphic ultrapowers.

Throughout the article, we will use the above notations. In addition, we will often consider the following subgroups of $T_0(\Gamma)$ or $T_1(\Gamma)$:
\[\widetilde\Gamma_n=\bigoplus_{i\geq n}\Gamma_i,\;\; \text{ and } \;\;\widetilde\Gamma_{n,n'}=\bigoplus_{n'>i\geq n}\Gamma_i,\;\;\text{ for every } n'>n\geq 1.\]

\subsection*{Acknowledgements}  I.C. is very grateful to Isaac Goldbring for kindly showing to him that the results in \cite{ZM69} can be used to produce a fourth elementary equivalence class of II$_1$ factors. A.I. would like to thank Sorin Popa and Thomas Sinclair for many stimulating discussions. We are also grateful to Ilijas Farah for various comments, and especially for pointing out Corollary \ref{maincor} to us.

\section{Preliminaries}

\subsection{Terminology} 
Throughout this article we work with {\it tracial von Neumann algebras} $(M,\tau)$, i.e. von Neumann algebras $M$ endowed with a faithful normal trace $\tau:M\rightarrow\mathbb C$.  We say that $M$ is {\it separable} if it is separable with respect to the norm $\|x\|_2=\tau(x^*x)^{1/2}$.
We denote by $\ms U(M)$ the group of {\it unitaries} of $M$.
If $n\geq 1$, then we denote by $M^{{\bar\otimes} n}$ the tensor product von Neumann algebra ${\bar\otimes}_{i=1}^nM$.
If $A,B\subset M$, then we denote $A'\cap B=\{x\in B\,|\,xy=yx,\;\text{for all}\;y\in A\}$.  A tracial von Neumann algebra $M$ is a {\it II$_1$ factor} if it is infinite dimensional and has trivial center.

If $\Gamma$ is a countable group, then we denote by $(u_g)_{g\in\Gamma}\subset\ms U(\ell^2(\Gamma))$ its left regular representation given by $u_g(\delta_h)=\delta_{gh}$, where $(\delta_h)_{h\in\Gamma}$ is the usual orthonormal basis of $\ell^2(\Gamma)$.
 The weak (operator) closure of the linear span of $(u_g)_{g\in\Gamma}$ is a tracial von Neumann algebra, which we denote by $L(\Gamma)$. The so-called {\it group von Neumann algebra} $L(\Gamma)$ is a II$_1$ factor precisely when $\Gamma$ has infinite non-trivial conjugacy classes (icc).

 \subsection{Ultrafilters and ultraproducts} 
In this subsection we collect together several elementary facts regarding ultraproducts of von Neumann algebras.
 
 An {\it ultrafilter} $\mathcal U$ on a set $S$ is collection of subsets  of $S$ which is closed under finite intersections, does not contain the empty set, and contains either $S'$ or $S\setminus S'$, for every subset $S'\subset S$. An ultrafilter $\mathcal U$ is called {\it free} if it contains the complements of all finite subsets of $S$.
 If $f\in\ell^{\infty}(S)$, then its limit along $\mathcal U$, denoted by $\lim_{\mathcal U}f(s)$, is the unique $\ell\in\mathbb C$ such that $\{s\in S\,|\, |f(s)-\ell|<\varepsilon\}\in\mathcal U$, for every $\varepsilon>0$. The map $\ell^{\infty}(S)\ni f\mapsto\lim_{\mathcal U}f(s)\in\mathbb C$ is a $*$-homomorphism, which allows to identify $\mathcal U$ with a point in the Stone-\v{C}ech compactification $\beta S$ of $S$. Via this identification, an ultrafilter $\mathcal U$ is free if and only if it belongs to $\beta S\setminus S$.

Given an ultrafilter $\cU$ on a set $S$ and a family of tracial von Neumann algebras $(M_s,\tau_s), s\in S$, we define the {\it ultraproduct} algebra $\prod_{\mathcal U}M_s$ as the quotient $\mathcal A/\mathcal I$, where $\mathcal A$ is the C$^*$-algebra $\mathcal A=\{(x_s)_s\in\prod {M_s}\,|\,\sup_s\|x_s\|<\infty\}$ and $\mathcal I$ is the closed ideal of $(x_s)_s\in\mathcal A$ such that $\lim_{\mathcal U}\|x_s\|_2=0$. It turns out that $\prod_{\mathcal U}M_s$ is a tracial von Neumann algebra, with the canonical trace given by $\tau((x_s)_s)=\lim_{\mathcal U}\tau_s(x_s)$. When $(M_s)_s$ is the constant family $(M)_s$, we write $M^\cU$ and call it the {\it ultrapower} von Neumann algebra. In this case, the map $\pi:M\rightarrow M^{\mathcal U}$ given by $\pi(x)=(x_s)_s$, where $x_s=x$, for all $s\in S$, is an injective $*$-homomorphism. 

Next, we recall the known fact that depending on the ultrafilter $\cU$, $M^{\cU}$ is either non-separable or isomorphic to $M$ (see Lemma \ref{complete}).

\begin{definition}
An ultrafilter $\mathcal U$ on a set $S$ is called {\it countably cofinal} (or {\it countably incomplete}) if there is a sequence $\{A_n\}$ in $\mathcal U$ such that $\cap_{n}A_n=\emptyset$. Otherwise, $\mathcal U$ is called {\it countably complete}. 
\end{definition}

Any free ultrafilter on a countable set is countably cofinal, while any principal (i.e.\ non-free) ultrafilter is countably complete.  
The hypothesis that there exists a countably complete free ultrafilter is a very strong axiom which is not provable from ZFC (e.g.\ \cite[Section 5]{Ke10}).
However, such set-theoretic issues will not be important here.

\begin{lemma}\label{ccofinal}  Let $\mathcal U$ be a countably cofinal ultrafilter on a set $S$. 
 For every $s\in S$, let $M_s$ be a tracial von Neumann algebra and $\{M_{s,n}\}_{n\geq 1}$ be an increasing sequence of von Neumann subalgebras whose union is weakly dense in $M_s$. Let $Q\subset \prod_{\mathcal U}M_s$ be a separable subalgebra. 

Then for every $s\in S$ we can find an integer $n_s\geq 1$ such that $Q\subset\prod_{\mathcal U}M_{s,n_s}$.
\end{lemma}

{\it Proof.} Since $\mathcal U$ is countably cofinal, we can find a sequence $\{A_n\}_{n\geq 2}$ in $\mathcal U$ such that $\cap_{n}A_n=\emptyset$. Let $A_1=S\setminus A_2$.
For $s\in S$, let $f(s)$ be the largest integer $n\geq 1$ such that $s\in A_n$.  It is clear that $f:S\rightarrow\mathbb N$ is well-defined and $\lim_{\mathcal U}f(s)=+\infty$.

Let $Q\subset\prod_{\mathcal U}M_s$ be a separable subalgebra. Let $\{x_k\}_{k\geq 1}$ be a $\|\cdot\|_2$-dense sequence in $Q$. 
For $k\geq 1$, represent $x_k=(x_{k,s})_s$, where $x_{k,s}\in M_s$, for all $s\in S$. Let $s\in S$. Since $\cup_{n\geq 1}M_{s,n}$ is weakly dense in $M$ and $M_{s,n}\subset M_{s,n+1}$, for all $n\geq 1$, we can find $n_s\geq 1$ such that $$\|x_{k,s}-E_{M_{s,n_s}}(x_{k,s})\|_2\leq\frac{1}{f(s)},\;\;\text{for all $1\leq k\leq f(s)$}.$$

Since $\lim_{\mathcal U}f(s)=+\infty$, it follows that for every $k\geq 1$ we have $\lim_{\mathcal U}\|x_{k,s}-E_{M_{s,n_s}}(x_{k,s})\|_2=0$. This implies that $x_k\in\prod_{\mathcal U}M_{s,n_s}$, for every $k\geq 1$, and hence $Q\subset\prod_{\mathcal U}M_{s,n_s}$. \hfill$\blacksquare$

The first assertion of the next lemma is well-known \cite{Fe56}, while the second assertion follows from the proof of \cite[Proposition 6.1(2)]{GH01}. Nevertheless, we include a proof for completeness. 

\begin{lemma}\label{complete}
Let $\mathcal U$ be an ultrafilter on a set $S$. Let $(M,\tau)$ be a tracial von Neumann algebra. 
\begin{enumerate}
\item If $\mathcal U$ is countably cofinal and $M$ has a diffuse direct summand, then $M^{\mathcal U}$ is non-separable.
\item If  $\mathcal U$ is countably complete and $M$ is separable, then $\pi:M\rightarrow M^{\mathcal U}$ given by $\pi(x)=(x_s)_s$, where $x_s=x$, for all $s\in S$, is a $*$-isomorphism.
\end{enumerate}
\end{lemma}

{\it Proof.} (1) If $p\in M$ is a projection, then $(pMp)^{\mathcal U}$ is a subalgebra of $M^{\mathcal U}$. If $M$ is a diffuse von Neumann algebra, then $M$ contains a copy of $A:=L^{\infty}([0,1])$, hence $M^{\mathcal U}$ contains a copy of $A^{\mathcal U}$. 

We may therefore reduce to the case when $M=A$. 
Let $\{A_n\}_{n\geq 1}$ be an increasing sequence of finite dimensional subalgebras of $A$ such that $\cup_n A_n$ is weakly dense in $A$. Assuming that $A^{\mathcal U}$ is separable, Lemma \ref{ccofinal} implies that $A^{\mathcal U}=\prod_{\mathcal U}A_{n_s}$, for some integers $n_s\geq 1, s\in S$.
But since $A$ is diffuse and $A_{n_s}$ is finite dimensional, we can find $u_s\in\ms U(A)$ such that $E_{A_{n_s}}(u_s)=0$, for every $s\in S$. Then $u=(u_s)_s\in\ms U(A^{\mathcal U})$ would be orthogonal to $\prod_{\mathcal U}A_{n_s}$, which is a contradiction.

\noindent (2) Since $\mathcal U$ is countably complete, $\cap_n A_n\not=\emptyset$, for any sequence $\{A_n\}$ in $\mathcal U$. Then the collection $\mathcal U'$ of all sets of the form $\cap_n A_n$, where $\{A_n\}$ is a sequence in $\mathcal U$, is a filter on $S$. Since $\mathcal U\subset\mathcal U'$ and  $\mathcal U$ is an ultrafilter, we get that $\mathcal U'=\mathcal U$ and hence $\cap_n A_n\in\mathcal U$, for any sequence $\{A_n\}$ in $\mathcal U$.

Let $f_m\in\ell^{\infty}(S)$, for $m\geq 1$, and put $\ell_m=\lim_{\mathcal U}f_m(s)$. Then the previous paragraph implies that \begin{equation}\label{inter}\{s\in S \,|\, f_m(s)=\ell_m \,|\;\text{ for all }m\geq 1\}=\bigcap_{m,n\geq 1}\{s\in S\,|\, |f_m(s)-\ell_m|<\frac{1}{n}\}\in\mathcal U.\end{equation}

Assuming that $M$ is separable, let us show that $\pi:M\rightarrow M^{\mathcal U}$ is onto. Indeed, let $x=(x_s)_s\in M^{\mathcal U}$.  Let $\{z_m\}$ be a $\|\cdot\|_2$-dense sequence in $M$. For every $m$, define $f_m\in\ell^{\infty}(S)$ by  $f_m(s)=\tau(x_sz_m)$. By \eqref{inter}, there exists a set $A\in\mathcal U$ such that $f_m(s)=f_m(s')$, for all $s,s'\in A$ and every $m\geq 1$. Hence $x_s=x_s'$, for every $s,s'\in A$. Choosing $s_0\in A$, it clearly follows that $\pi(x_{s_0})=x$. This shows that $\pi$ is onto. Since $\pi$ is also injective, we conclude that $\pi$ is a $*$-isomorphism.
\hfill$\blacksquare$

Let us also record a simple consequence of Lemma \ref{ccofinal}, more specific to our problem. 

\begin{corollary}\label{combo2}
Let $\mathcal U$ be a countably cofinal ultrafilter on a set $S$. Let $\Gamma$ be a countable group, $\alpha\in\{0,1\}$, and consider the notation from Section \ref{sectionconstruction}. Denote by $M=L(T_{\alpha}(\Gamma))$ and $P_n=L(\widetilde\Gamma_n)$, for $n\geq 1$. For every $s\in S$, let $t_s\geq 1$ be an integer and let $Q_s$ be a tracial von Neumann algebra. 

If $A\subset\prod_{\cU}(M^{{\bar\otimes} t_s}{\bar\otimes} Q_s)$ is a separable subalgebra, then there are integers $n_s\geq 1$, $s\in S$, satisfying \[\prod_{\cU}P_{n_s}^{{\bar\otimes} t_s}\subset A'\cap\prod_{\cU} M^{{\bar\otimes} t_s}.\]
\end{corollary}

{\it Proof.}  From the definition of $T_0$ and $T_1$, we see that the increasing union $\cup_{n\geq 1}(P_n'\cap M)$ is weakly dense in $M$.  Thus, the increasing union $\cup_{n\geq 1}\big[\big((P_n^{{\bar\otimes} t_s})'\cap M^{{\bar\otimes} t_s}\big){\bar\otimes} Q_s\big]$ is  weakly dense in $M^{{\bar\otimes} t_s}{\bar\otimes} Q_s$, for all $s\in S$.   
Using Lemma \ref{ccofinal}, for each $s\in S$, there exists an integer $n_s\geq 1$ such that $A\subset\prod_{\mathcal U}\big[\big((P^{{\bar\otimes} t_s}_{n_s})'\cap M^{{\bar\otimes} t_s}\big){\bar\otimes} Q_s\big].$
This clearly implies the conclusion.
\hfill$\blacksquare$

\subsection{Residual inclusions} A subalgebra $P$ of a tracial von Neumann algebra $M$ is called residual if it ``absorbs" central sequences: any central sequence of $M$ asymptotically lies in $P$. In this subsection, we define and use a quantitative notion of residual subalgebras.

\begin{definition} Let $(M,\tau)$ be a tracial von Neumann algebra, $k\geq 1$ an integer, and $C>0$.
A von Neumann subalgebra $P\subset M$ is said to be $(k,C)$-{\it residual} if  there exist unitary elements $u_1, u_2, ..., u_k\in M$ such that for all $\xi\in M$ we have $$
\|\xi-E_{P}(\xi)\|_2^2\leq C \sum^k_{i=1}\|[\xi, u_i]\|_2^2.$$
\end{definition}

\begin{lemma} \label{bigger}
Let  $k\geq 1$ and $C>0$.   Let $\mathcal U$ be an ultrafilter on a set $S$. For any $s\in S$, let $M_s$ and $Q_s$ be tracial von Neumann algebras, and $P_s\subset M_s$ be a  $(k,C)$-residual von Neumann subalgebra.

Then there exists a separable von Neumann subalgebra  $A\subset \prod_{\mathcal U} M_s$ such that $$A'\cap\prod_{\mathcal U} (M_s{\bar\otimes} Q_s)\subset  \prod_{\mathcal U}(P_s{\bar\otimes} Q_s).$$
\end{lemma}

{\it Proof.} Let $s\in S$. Let $u^{s}_1,... ,u^{s}_k\in\ms U(M_s)$ such that
$\|\xi-E_{P_{s}}(\xi)\|_2^2\leq C \sum^k_{i=1}\|[\xi,u_i^s]\|_2^2$,  for all $\xi\in M_s$. From this it follows that $\|\xi-E_{P_s{\bar\otimes} Q_s}(\xi)\|_2^2\leq C\sum_{i=1}^k\|[\xi,u^s_i]\|_2^2$, for all $\xi\in M_s{\bar\otimes} Q_s$.
Denote by $u_i=(u^{s}_i)_s\in \ms U(\prod_{\mathcal U} M_s)$, for $1\leq i\leq k$. Then the last inequality implies that $$\|\xi-E_{\prod_{\mathcal U}(P_{s}{\bar\otimes} Q_s)}(\xi)\|_2^2\leq C \sum^k_{i=1}\| [\xi,u_i]\|_2^2,\;\;\text{for all $\xi\in\displaystyle{\prod_{\mathcal U}}(M_s{\bar\otimes} Q_s)$}.$$ 

Finally, we notice the von Neumann algebra $A$ generated by $u_1,...,u_k$ satisfies the conclusion. \hfill$\blacksquare$

\begin{definition}\cite{MD69a}\label{sres} Let $\Gamma$ be a countable group. A subgroup $\Lambda<\Gamma$ is called \emph{strongly residual} if there exist elements $a,b\in \Gamma$ and a subset  $F\subset\Gamma\setminus\Lambda$ satisfying the following properties:     
\begin{enumerate}
\item [(i)] $a\Lambda a^{-1}=\Lambda$,
\item [(ii)] $aFa^{-1}\cup F=\Gamma\setminus\Lambda$, and
\item [(iii)] $\{ b^{k}F b^{-k}\}_{k\in \mathbb Z}$ is a family of disjoint subsets of $\Gamma\setminus \Lambda$.
\end{enumerate}
\end{definition}

\begin{lemma}\label{resalg}
Let $\Lambda<\Gamma$ be a strongly residual subgroup.

Then $L(\Lambda)\subset L(\Gamma)$ is a  $(2,100)$-residual subalgebra.
\end{lemma}

{\it Proof.} 
Let $a,b\in\Gamma$ and $F\subset\Gamma\setminus\Lambda$ be as in Definition \ref{sres}.
Let  $\xi\in L(\Gamma)$. View $\xi\in\ell^2(\Gamma)$ and for any subset $A\subset\Gamma$, define $\nu(A)=\sum_{h\in A}|\xi(h)|^2$. Then $\nu$ is a finite measure on $\Gamma$ and the Cauchy-Schwarz inequality implies that for every $g\in\Gamma$ and $A\subset\Gamma$ we have
\begin{align*}
|\nu(gAg^{-1})-\nu(A)| & = \Big|\sum_{h\in A}|\xi(ghg^{-1})|^2-|\xi(h)|^2\Big| \\
& \leq \sum_{h\in A} (|\xi(ghg^{-1})| + |\xi(h)|)||\xi(ghg^{-1})| - |\xi(h)||\\
& \leq (\nu(gAg^{-1})^{1/2}+\nu(A)^{1/2})\sum_{h \in A} |\xi(ghg^{-1}) - \xi(h)|^2. 
\end{align*}
Hence  
 \begin{equation}\label{ecu1}|\nu(gAg^{-1})-\nu(A)| \leq (\nu(gAg^{-1})^{1/2}+\nu(A)^{1/2})\|[u_g,\xi]\|_2.\end{equation}

On the other hand, by using conditions (ii) and (iii) from Definition \ref{sres} we get that 
\begin{equation}\label{ecu2}\nu(F)+\nu(aFa^{-1})\geq\nu(\Gamma\setminus\Lambda)\geq\nu(F)+\nu(bFb^{-1})+\nu(b^{-1}Fb).\end{equation}
Combining \eqref{ecu1} with \eqref{ecu2}, and using $F,aFa^{-1},bFb^{-1},b^{-1}Fb\subset\Gamma\setminus\Lambda$ we deduce that \begin{align*}\nu(\Gamma\setminus\Lambda)&\leq  3(\nu(F)+\nu(aFa^{-1}))-2(\nu(F)+\nu(bFb^{-1})+\nu(b^{-1}Fb))
\\ & \leq 3|\nu(aFa^{-1})-\nu(F)|+2|\nu(bFb^{-1})-\nu(F)|+2|\nu(b^{-1}Fb)-\nu(F)|\\&\leq 6\;\nu(\Gamma\setminus\Lambda)^{1/2}\|[u_a,\xi]\|_2+8\;\nu(\Gamma\setminus\Lambda)^{1/2}\|[u_b,\xi]\|_2\\&\leq 10\;\nu(\Gamma\setminus\Lambda)^{1/2}(\|[u_a,\xi]\|_2^2+\|[u_b,\xi]\|_2^2)^{1/2}.\end{align*}

Since $\nu(\Gamma\setminus\Lambda)=\|\xi-E_{L(\Lambda)}(\xi)\|_2^2$, the conclusion follows.
\hfill$\blacksquare$

\begin{lemma}\label{SRSex}
Let $\Gamma$ be a countable group,  and consider the notations from Section \ref{sectionconstruction}. 

Then $\widetilde\Gamma_n$ is strongly residual in $T_{\alpha}(\Gamma)$, for every $n\geq 1$ and $\alpha\in\{0,1\}$.

\end{lemma}

This statement is proven in \cite[\S 21]{DL69} in the case $\alpha=0$ and $n=1$, and is used (without proof) in full generality in \cite{MD69a,MD69b}. For completeness, we provide a proof.

{\it Proof.} Let $n\geq 1$. When $\alpha=0$, we define $\Sigma_n<T_0(\Gamma)$ to be the subgroup generated by $\widetilde\Gamma,\Lambda_1,\Lambda_2,...,\Lambda_n$, and $\Delta_n<T_0(\Gamma)$ to be the subgroup generated by $\widetilde\Gamma_n,\Lambda_{n+1},\Lambda_{n+2},...$.  
Similarly, when $\alpha=1$, we define $\Sigma_n<T_1(\Gamma)$ to be the subgroup generated by $\widetilde\Gamma\rtimes S_{\infty},\Lambda_1,\Lambda_2,...,\Lambda_n$, and $\Delta_n<T_1(\Gamma)$ to be the subgroup generated by $\widetilde\Gamma_n,\Lambda_{n+1},\Lambda_{n+2},...$. 
In both cases one can check that $\Sigma_n$ and $\Delta_n$ generate $T_{\alpha}(\Gamma)$ so that $T_{\alpha}(\Gamma)=\Sigma_n \ast_{\widetilde\Gamma_n}\Delta_n$. Moreover, if $a$ and $b$ are the generators of $\Lambda_1$ and $\Lambda_{n+1}$, then $a$ commutes with $\widetilde\Gamma_n$, $a\in\Sigma_n\setminus\widetilde\Gamma_n$, and $b^k\in\Delta_n\setminus\widetilde\Gamma_n$, for all $k\geq 1$.

The conclusion is now a consequence of the following fact. Let $G = H_1 \ast_{H_0} H_2$ be an amalgamated free product group such that  there exist $a\in H_1\setminus H_0$, $b\in H_2\setminus H_0$ satisfying $aH_0a^{-1}= H_0$ and $b^k\notin H_0$, for all $k\geq 1$. Then $H_0$ is a strongly residual subgroup of $G$. Indeed, let $F \subset G \setminus H_0$ be the set of reduced words of the form $g_1g_2...g_k$, for some $k\geq 1$ and $g_1 \in H_1\setminus H_0$, $g_2\in H_2\setminus H_0$, $g_3 \in H_1 \setminus H_0,...$. It is easy to see that $a,b,F$ verify conditions (i)-(iii) listed in Definition \ref{sres}.
\hfill$\blacksquare$

\begin{lemma}\emph{\cite{MD69a}}\label{srprod} 
 Let $\Lambda_i<\Gamma_i$ be strongly residual subgroups, for every $1\leq i\leq n$.
  
 Then $\oplus^n_{i=1} \Lambda_i< \oplus^n_{i=1} \Gamma_i$ is a strongly residual subgroup.
\end{lemma}

For completeness, we recall the proof from \cite{MD69a}. 

{\it Proof.}  Let $a_i,b_i\in\Gamma_i$ and $F_i\subset\Gamma_i\setminus\Lambda_i$ be as in Definition \ref{sres}. Denote by $\Lambda=\oplus^n_{i=1} \Lambda_i$ and $\Gamma=\oplus^n_{i=1} \Gamma_i$. Let $a=(a_1,...,a_n)$, $b=(b_1,...,b_n)\in\Gamma$, and consider the following set 
\[F := \{(g_1,...,g_n) \in \Gamma \, | \,  \exists i \text{ such that } g_1\in\Lambda_1, \dots, g_{i-1}\in\Lambda_{i-1} \text{ and } g_i \in F_i \}.\] 
Then  $a,b,F$ satisfy conditions (i)-(iii) from Definition \ref{sres} for $\Lambda<\Gamma$. \hfill$\blacksquare$

The following key corollary of the above results will be frequently used in the sequel.

\begin{corollary}\label{combo}
Let $\mathcal U$ be an ultrafilter on a set $S$. Let $\Gamma$ be a countable group, $\alpha\in\{0,1\}$, and consider the notation from Section \ref{sectionconstruction}. Denote by $M=L(T_{\alpha}(\Gamma))$ and $P_n=L(\widetilde\Gamma_n)$, for every $n\geq 1$. For every $s\in S$, let $n_s,t_s\geq 1$ be integers and let $Q_s$ be a tracial von Neumann algebra. 

Then there exists a separable subalgebra $A\subset\prod_{\cU}M^{{\bar\otimes} t_s}$ such that $$A'\cap\prod_{\cU}(M^{{\bar\otimes} t_s}{\bar\otimes} Q_s)\subset\prod_{\cU}(P_{n_s}^{{\bar\otimes} t_s}{\bar\otimes} Q_s).$$
\end{corollary}

{\it Proof.} Let $s\in S$.
By combining Lemmas \ref{SRSex} and \ref{srprod}, we get that $\oplus_{i=1}^{t_s}\widetilde\Gamma_{n_s}<\oplus_{i=1}^{t_s}T_{\alpha}(\Gamma)$ is a strongly residual subgroup. Then Lemma \ref{resalg} implies that $P_{n_s}^{{\bar\otimes} t_s}\subset M^{{\bar\otimes} t_s}$ is a $(2,100)$-residual subalgebra. The conclusion now follows from Lemma \ref{bigger}. 
\hfill$\blacksquare$

\section{McDuff depth of a von Neumann algebra}
\label{sectiondepth}

\subsection{Properties at depth $k$}

\begin{definition}
Let $\cM$ be a (typically non-separable) von Neumann algebra, and $I$ be a directed set. On the set of subalgebras of $\cM$, consider the partial order given by inclusion. A decreasing net $(A_i)_{i \in I}$ of subalgebras of $\mathcal M$ is called a {\it residual net} if:
\begin{itemize}
\item For any separable subalgebra $Q \subset \cM$, there exists $i \in I$ such that $A_i \subset Q' \cap \cM$, and
\item For any $i \in I$, there exists a separable subalgebra $Q \subset \cM$ such that $Q' \cap \cM \subset A_i$.
\end{itemize}
A residual net is called {\it trivial} if there exists $i \in I$ such that $A_i = \C$. 
\end{definition}

\begin{example}\label{standard}
Given a von Neumann algebra $\cM$, consider the set $I$ of separable subalgebras of $\cM$, ordered by inclusion. Then the net $(Q' \cap \cM)_{Q \in I}$ is clearly a residual net, which we call the {\it standard residual net}.
\end{example}

The following result provides our main example of a residual net.

\begin{lemma}\label{MDex}
Let $\Gamma$ be a countable group and $\alpha \in \{0,1\}$. For $n'>n\geq 1$, let $M = L(T_\alpha(\Gamma))$, $P_n = L(\widetilde\Gamma_n)$, $P_{n,n'}=L(\widetilde\Gamma_{n,n'})$, where $\widetilde\Gamma_{n,n'}<\widetilde\Gamma_n<\widetilde\Gamma<T_{\alpha}(\Gamma)$ are defined as in Section \ref{sectionconstruction}. 

Let $\mathcal U$ be a countably cofinal ultrafilter on a set $S$. For every $s\in S$, let $t_s\geq 1$ be an integer. Define $\cM=\prod_{\mathcal U}M^{{\bar\otimes} t_s}$.
Endow $I = \N^{S}$ with the following partial order:
$(n_s)_s < (m_s)_s \, \text{ iff } \, n_s < m_s, \, \text{ for all } \, s\in S.$ 
For $i =(n_s)_{s\in S}\in I$, define $A_i=\prod_\cU {P_{n_s}}^{{\bar\otimes} t_s}$.

Then $(A_i)_{i \in I}$ is a residual net of $\cM$. Moreover, if $\Gamma$ is icc and $i=(n_s)_s<j=(m_s)_s$, then  $A_j' \cap A_i=\prod_\cU {P_{n_s,m_s}}^{{\bar\otimes} t_s}$.
\end{lemma}
{\it Proof.}
The fact that $(A_i)_{i\in I}$ is a residual net of $\cM$ follows easily from Lemmas \ref{combo2} and \ref{combo}.
Now, if $\Gamma$ is icc, then $L(\Gamma)$ is a II$_1$ factor. It follows that for $n'>n\geq 1$ we have $P_{n'}'\cap P_n=P_{n,n'},$ which clearly implies the moreover part.
\hfill$\blacksquare$

Motivated by the moreover assertion of Lemma \ref{MDex}, we introduce the following definition.

\begin{definition}
Let $\cP$ be a property for von Neumann algebras.
Let $\cM$ be a von Neumann algebra with a residual net $(A_i)_{i\in I}$. We say that $\cM$ has {\it property $\cP$ at depth $1$} if for all $i_1 \in I$, there exists $i_2 > i_1$ such that for all $i_3 > i_2$ there exists $i_4 > i_3$, such that the inclusion $A_{i_3}' \cap A_{i_2} \subset A_{i_4}' \cap A_{i_1}$ contains an intermediate von Neumann subalgebra with property $\cP$.
\end{definition}

\begin{remark}\label{trivial}
A von Neumann algebra $\cM$ is trivial at depth $1$ if and only if it admits a separable subalgebra $Q\subset\cM$ such that $Q'\cap\cM=\mathbb C1$.
\end{remark}

\begin{lemma}\label{indepnet}
Having property $\cP$ at depth $1$ is independent of the choice of a residual net.
\end{lemma}
{\it Proof.}
Consider a von Neumann algebra $\cM$ with two residual nets $(A_i)_{i \in I}$ and $(B_j)_{j \in J}$. From the definition of residual nets we see that for every $i \in I$ there exists $j \in J$ such that $B_j \subset A_i$. Symmetrically, for every $j\in J$ there exists $i\in I$ such that $A_i\subset B_j$.

Assume that $\cM$ has property $\cP$ at depth $1$ with respect to $(A_i)_{i\in I}$.

Fix $j_1 \in J$. Then there exists $i_1 \in I$ such that $A_{i_1} \subset B_{j_1}$. Take $i_2 > i_1$ as in the definition of property $\cP$ at depth $1$.
Then there exists $j_2 > j_1$ such that $B_{j_2} \subset A_{i_2}$.

Take an arbitrary $j_3 > j_2$, and pick $i_3>i_2$ such that $A_{i_3} \subset B_{j_3}$. Next, we find $i_4 > i_3$ as in the definition of property $\cP$ at depth $1$.
Then there exists $j_4 > j_3$ such that $B_{j_4} \subset A_{i_4}$.

Altogether, we have the following inclusions
\[B_{j_4} \subset A_{i_4} \subset A_{i_3} \subset B_{j_3} \subset B_{j_2} \subset A_{i_2} \subset A_{i_1} \subset B_{j_1}.\]
From this we get that
\[B_{j_3}' \cap B_{j_2} \subset A_{i_3}' \cap A_{i_2} \subset A_{i_4}' \cap A_{i_1} \subset B_{j_4}' \cap B_{j_1}.\]
By our choice of $i_1,i_2,i_3,i_4 \in I$, there is an intermediate subalgebra with property $\cP$ inside the inclusion 
$A_{i_3}' \cap A_{i_2} \subset A_{i_4}' \cap A_{i_1}$. Thus, $\cM$ has property $\cP$ at depth $1$ with respect to  $(B_j)_{j\in J}$.
\hfill$\blacksquare$

\begin{definition}\label{atdepthk}
Let $\cP$ be a property of von Neumann algebras. We define inductively on $k \geq 0$ what it means for a von Neumann algebra $\cM$ to have {\it property $\cP$ at depth $k$}. We denote this property by $\cP^{(k)}$.
\begin{itemize}
\item If $k = 0$, then we say that $\cM$ has property $\cP^{(0)}$  if it has property $\cP$.
\item If $k \geq 0$, then we say that $\cM$ has $\cP^{(k+1)}$  if it has property $\cP^{(k)}$ at depth $1$.
\end{itemize}
\end{definition}

\begin{definition}
A von Neumann algebra $\cM$ is said to have {\it finite McDuff depth} if there exists $k$ such that $\cM$ is trivial at depth $k$. The {\it McDuff depth} of $\cM$ is defined as the smallest $k \geq 0$ such that $\cM$ is trivial at depth $k+1$.
If $\cM$ does not have finite McDuff depth, then we define its McDuff depth to be infinite.
\end{definition}

\begin{examples} Let $M$ be a separable II$_1$ factor and $\cU$ be a countably cofinal ultrafilter.
\begin{itemize}
\item By Remark \ref{trivial}, $M^{\cU}$ has depth $0$ if and only if $M$ does not have property Gamma.
\item If $M$ has property Gamma but is non McDuff, then $M' \cap M^\cU$ is abelian and non trivial. This easily implies that $M^\cU$ has infinite depth.
\item If $M$ is the hyperfinite II$_1$ factor, then $M^\cU$ has depth $1$.
\end{itemize}
\end{examples}

\subsection{Computing the depth} The aim of this subsection is to prove the following result.

\begin{theorem}\label{countable} 
Let $\cU$ be a countably cofinal ultrafilter, and $\Gamma$ be a non-trivial countable group. Let $\pmb\alpha$ be a (finite or infinite) sequence of $0$'s and $1$'s. Let $M_{\pmb\alpha}(\Gamma)$ be as defined in Section \ref{sectionconstruction}.

Then the depth of $M_{\pmb\alpha}(\Gamma)^\cU$ is at least the length of $\pmb\alpha$. Moreover, if $\Gamma =\mathbb F_2$, then we have equality.
\end{theorem}

Let us fix a countably cofinal ultrafilter $\cU$ on a set $S$. Towards proving Theorem \ref{countable}, we first provide an upper bound on depth when $\Gamma = \mathbb F_2$.

\begin{lemma}\label{count1}
For all $k \geq 0$, $\pmb\alpha \in \{0,1\}^k$, and all integers $t_s \geq 1, s\in S$, the von Neumann algebra $\prod_{\cU} M_{\pmb\alpha}(\mathbb F_2)^{{\bar\otimes} t_s}$ is trivial at depth $k+1$.
\end{lemma}
{\it Proof.}
We proceed by induction on the length $k$ of the sequence $\pmb\alpha$. First assume $k = 0$. Then $\pmb\alpha$ is the empty sequence and $M_{\pmb\alpha}(\mathbb F_2) = L(\mathbb F_2)$. Since the trivial subgroup is strongly residual in $\mathbb F_2$ (see e.g.\ the proof of Lemma \ref{SRSex}), by combining Lemmas \ref{srprod}, \ref{resalg}, and \ref{bigger}, we deduce the existence of a separable subalgebra $Q\subset\prod_{\cU} M_{\pmb\alpha}(\mathbb F_2)^{{\bar\otimes} t_s}$ with trivial relative commutant. In other words, $\prod_{\cU} M_{\pmb\alpha}(\mathbb F_2)^{{\bar\otimes} t_s}$ is trivial at depth $1$ (see Remark \ref{trivial}).

Assume the conclusion holds for $k \geq 0$. Take integers $t_s\geq 1, s\in S$ and a sequence $\pmb\alpha \in \{0,1\}^{k+1}$ of length $k+1$. Put $\cM = \prod_{\cU} M_{\pmb\alpha}(\mathbb F_2)^{{\bar\otimes} t_s}$.
We want to show that, at depth 1, $\cM$ has the property of being trivial at depth $k+1$.  

By Lemma \ref{indepnet}, in order to check this property at depth 1 we can use any residual net for $\cM$. 
From Definition \ref{defintro} we see that $T_{\pmb\alpha}(\mathbb F_2)=T_{\alpha_1}(T_{\pmb\beta}(\mathbb F_2))$ and hence $M_{\pmb\alpha}(\mathbb F_2) = L(T_{\alpha_1}(T_{\pmb\beta}(\mathbb F_2)))$, where $\pmb\beta = (\alpha_2,\dots,\alpha_{k+1}) \in \{0,1\}^k$.
Applying Lemma \ref{MDex} (to $\Gamma = T_{\pmb\beta}(\mathbb F_2)$ and $\alpha = \alpha_1$) we obtain  that
 $A_i = \prod_\cU {P_{n_s}}^{{\bar\otimes} t_s}$, where $i = (n_s)_s\in I$, is a residual net for $\cM$.

Now, for any indices $i_4 > i_3 > i_2 > i_1$, the inclusion $A_{i_3}' \cap A_{i_2} \subset A_{i_4}' \cap A_{i_1}$ has $A_{i_4}' \cap A_{i_1} $ as an intermediate algebra. Moreover, since $T_{\pmb\beta}(\mathbb F_2))$ is icc, by Lemma \ref{MDex} we get that $A_{i_4}'\cap A_{i_1}=\prod_\cU {P_{{n_1}_s,{n_4}_s}}^{{\bar\otimes} t_s}$. Since $P_{m,m'}\cong M_{\pmb\beta}(\mathbb F_2)^{{\bar\otimes} (m'-m)}$, for any $m'>m\geq 1$, we conclude that $A_{i_4}'\cap A_{i_1}$  is of the form $\prod_\cU {M_{\pmb\beta}(\mathbb F_2)}^{{\bar\otimes} v_s}$, for some integers $v_s\geq 1, s\in S$. By our induction assumption, $A_{i_4}'\cap A_{i_1}$ is trivial at depth $k+1$. This shows that $\cM$ is trivial at depth $k+2$, as desired.
\hfill$\blacksquare$

Obtaining a lower bound on depth requires additional work. We start by defining some more general residual nets.

\begin{definition}
Let $\cN \subset \cM$ be an inclusion of von Neumann algebras and let $I$ be a directed set. Let $(A_i)_{i \in I}$ and $(B_i)_{i \in I}$ be two decreasing nets of subalgebras of $\cM$ such that $A_i \subset \cN$, $B_i \subset \cM$, and $A_i \subset B_i$, for all $i \in I$. The net $(A_i \subset B_i)_{i\in I}$ is called a {\it residual pair} for the inclusion $\cN\subset \cM$ if the following two properties are satisfied:
\begin{itemize}
\item For any $i \in I$, there exists a separable subalgebra $Q \subset \cN$ such that $Q' \cap \cM \subset B_i$, and
\item For any separable subalgebra $Q \subset \cM$, there exists $i \in I$ such that $A_i \subset Q' \cap \cN$.
\end{itemize}
\end{definition}

If $\cN = \cM$ and $(A_i)_{i\in I}$ is a residual net of $\cM$, then $(A_i \subset A_i)_{i\in I}$ is a residual pair for $\cN\subset\cM$. 
As in Lemma \ref{MDex}, we have the following key example.

\begin{lemma}\label{MDex2}
Let $\Gamma$ be a countable group and $\alpha \in \{0,1\}$. For $n'>n\geq 1$, let $M = L(T_\alpha(\Gamma))$, $P_n = L(\widetilde\Gamma_n)$, $P_{n,n'}=L(\widetilde\Gamma_{n,n'})$, where $\widetilde\Gamma_{n,n'} < \widetilde\Gamma_n < \widetilde\Gamma < T_{\alpha}(\Gamma)$ are defined as in Section \ref{sectionconstruction}. 

Recall that $\mathcal U$ is a countably cofinal ultrafilter on the set $S$. Let $Q_s, s\in S$, be tracial von Neumann algebras. Define $\cN=\prod_{\mathcal U}M$ and $\cM=\prod_{\mathcal U}(M{\bar\otimes} Q_s)$.
Consider $I = \N^{S}$ ordered as in Lemma \ref{MDex}.
For $i=(n_s)_{s\in S}\in I$, define $A_i=\prod_\cU {P_{i_s}} \subset \cN$ and $B_i=\prod_\cU (P_{i_s}{\bar\otimes} Q_s) \subset \cM$.

Then the net $(A_i\subset B_i)_{i \in I}$ is a residual pair for the inclusion $\cN\subset\cM$.

Moreover, if $\Gamma$ is icc and $i=(n_s)_s < j=(m_s)_s$, then we have  $A_j' \cap B_i=\prod_\cU (P_{n_s,m_s}{\bar\otimes} Q_s)$ and $B_j'\cap A_i=\prod_\cU P_{n_s,m_s}$.
\end{lemma}

The following result is a simple variation of Lemma \ref{indepnet}.

\begin{lemma}\label{inclusion}
Let $\cP$ be a property of von Neumann algebras.
Let $\cA \subset \cM \subset \cB$ be von Neumann algebras. Assume that $\cM$ has property $\cP$ at depth $1$. 

Then for any residual pair $(A_i \subset B_i)_{i \in I}$ for the inclusion $\cA \subset \cB$, we have:

For all $i_1 \in I$, there exists $i_2 > i_1$ such that for all $i_3 > i_2$ there exists $i_4 > i_3$ such that the inclusion $B_{i_3}' \cap A_{i_2} \subset A_{i_4}' \cap B_{i_1}$ contains an intermediate von Neumann algebra with property $\cP$.
\end{lemma}
{\it Proof.} 
Fix $i_1 \in I$. Then there exists a separable subalgebra $Q_1 \subset \cA$ such that $Q_1' \cap \cB \subset B_{i_1}$.

Since $\mathcal M$ has property $\cP$ at depth $1$,  Lemma \ref{indepnet} implies that $\cM$ has property $\cP$ at depth $1$ with respect to the standard residual net from Example \ref{standard}. Thus, there exists a separable subalgebra $\cM\supset Q_2\supset Q_1$ such that for any separable subalgebra $\cM\supset Q_3\supset Q_2$ we can find a separable subalgebra $\cM\supset Q_4\supset Q_3$ such that the inclusion $(Q_3'\cap \cM)' \cap (Q_2' \cap \cM) \subset (Q_4' \cap \cM)' \cap (Q_1' \cap \cM)$ contains an intermediate subalgebra with property $\cP$.

Next, let $i_2 > i_1$ such that $A_{i_2} \subset Q_2' \cap \cA$. 
Pick any index $i_3\in I$ with $i_3 > i_2$. Then one can find a separable subalgebra $Q \subset \cA$ such that $Q' \cap \cB \subset B_{i_3}$. 
Let $Q_3$ be the von Neumann algebra generated by $Q$ and $Q_2$. Then $\cM\supset Q_3\supset Q_2$ and $Q_3'\cap\cB\subset B_{i_3}$. 
Let $Q_4\supset Q_3$ be as given by the previous paragraph.  Since $Q_4\subset\mathcal M\subset\mathcal B$, there exists $i_4 > i_3$ such that $A_{i_4} \subset Q_4' \cap \cA$. 

Altogether, we have the following inclusions
\[A_{i_4} \subset Q_4' \cap \cM \subset Q_3' \cap \cM \subset B_{i_3}\]
\[A_{i_2} \subset Q_2' \cap \cM \subset Q_1' \cap \cM \subset B_{i_1}.\]
These lead to
\[B_{i_3}' \cap A_{i_2} \subset (Q_3'\cap \cM)' \cap (Q_2' \cap \cM) \subset (Q_4' \cap \cM)' \cap (Q_1' \cap \cM) \subset A_{i_4}' \cap B_{i_1}.\]

It follows that  the inclusion $B_{i_3}' \cap A_{i_2} \subset A_{i_4}' \cap B_{i_1}$ contains an intermediate von Neumann algebra with property $\cP$, as claimed.
\hfill$\blacksquare$

We are now ready to prove the second half of Theorem \ref{countable}: the lower bound on the depth.

\begin{lemma}\label{count2}
Fix a non-trivial countable group $\Gamma$.
Then for any $k\geq 0$ and $\pmb\alpha \in \{0,1\}^k$ and any family of tracial von Neumann algebras $Q_s$, $s\in S$, any intermediate von Neumann subalgebra $M_{\pmb\alpha}(\Gamma)^{\cU}=\prod_\cU M_{\pmb\alpha}(\Gamma) \subset \cM \subset \prod_\cU (M_{\pmb\alpha}(\Gamma) {\bar\otimes} Q_s)$ is not trivial at depth $k$.
\end{lemma}
{\it Proof.}
We proceed by induction on $k$. First assume $k = 0$. Then $\pmb\alpha$ is the empty sequence and $M_{\pmb\alpha}(\Gamma) = L(\Gamma) \neq \C1$, hence the ultrapower ${M_{\pmb\alpha}(\Gamma)}^\cU$ is not trivial, so $\cM$ is not trivial either.

Assume the result holds for some $k \geq 0$.  Take a sequence $\pmb\alpha$ of length $k + 1$ and suppose by contradiction there exist tracial von Neumann algebras $Q_s, s\in S$ and an intermediate von Neumann subalgebra $\cM$ satisfying  
\[M_{\pmb\alpha}(\Gamma)^\cU \subset \cM \subset \prod_\cU (M_{\pmb\alpha}(\Gamma) \bar\otimes Q_s),\]

that is trivial at depth $k+1$.

Therefore, at depth $1$, $\cM$ has the property of being trivial at depth $k$. From Definition \ref{defintro}, we see that $M_{\pmb\alpha}(\Gamma) =L(T_{\alpha_1}(T_{\pmb\beta}(\Gamma)))$, where $\pmb\beta = (\alpha_2,\dots,\alpha_{k+1}) \in \{0,1\}^k$. Let $(A_i \subset B_i)_{i \in I}$ be the residual pair for the inclusion $M_{\pmb\alpha}(\Gamma)^\cU \subset \prod_\cU (M_{\pmb\alpha}(\Gamma) \bar\otimes Q_s)$ obtained by applying Lemma \ref{MDex2} to $T_{\pmb\beta}(\Gamma)$ and $\alpha_1$ instead of $\Gamma$ and $\alpha$.

Lemma \ref{inclusion} then implies that for all $i_1 \in I$, there exists $i_2 > i_1$ such that for all $i_3 > i_2$ there exists $i_4 > i_3$ such that the inclusion ${B_{i_3}}' \cap A_{i_2} \subset {A_{i_4}}' \cap B_{i_1}$ contains an intermediate subalgebra which is trivial at depth $n$.
Since $T_{\pmb\beta}(\Gamma)$ is icc, Lemma \ref{MDex2} implies that for all indices $i_4 > i_3 > i_2 > i_1$ we have 
\[{B_{i_3}}' \cap A_{i_2} = \prod_\cU P_{n_{2,s},n_{3,s}} \, \text{ and } \, {A_{i_4}}' \cap B_{i_1} = \prod_\cU (P_{n_{1,s},n_{4,s}} \bar\otimes Q_s).\]

Choose $i_3=(n_{3,s})$ such that $n_{3,s} = n_{2,s}+1$ for all $s\in S$. Since $P_{m,m'} \cong {M_{\pmb\beta}(\Gamma)}^{\bar{\otimes} (m'-m)}$, for any $m'>m\geq 1$, we see that the inclusion ${B_{i_3}}' \cap A_{i_2} \subset {A_{i_4}}' \cap B_{i_1}$ is of the form
\[\prod_\cU P_{n_{2,s},n_{2,s}+1}= \prod_\cU M_{\pmb\beta}(\Gamma) \subset \prod_\cU (P_{n_{1,s},n_{4,s}} {\bar{\otimes}} Q_s) = \prod_\cU (M_{\pmb\beta}(\Gamma) {\bar{\otimes}} \widetilde Q_s),\]
for some tracial von Neumann algebras $\widetilde Q_s$. Since by the induction assumption there is no intermediate subalgebra in this inclusion which is trivial at depth $k$, we get a contradiction.
\hfill$\blacksquare$

\section{Distinguishing uncountably many ultrapowers}
\label{sectionpropV}

\subsection{Property $\wV$ and proof of the main results}

In order to show that the II$_1$ factors $M_{\pmb{\alpha}}(\Gamma)$ are non-isomorphic, McDuff introduced  \cite{MD69b} a certain property for separable II$_1$ factors, called property $V$ (cf.\ with the earlier notions of asymptotically abelian II$_1$ factors \cite{Sa68, DL69, ZM69}).  In this section, inspired by property $V$, we define the following new property for non-separable von Neumann algebras: 

\begin{definition}\label{propV}
A non-separable von Neumann algebra $\cM$ has {\it property $\wV$} if there exists a separable subalgebra $A \subset \cM$ such that for any separable subalgebra $B \subset A' \cap \cM$ and any separable subalgebra $C \subset \cM$,  there exists a unitary $u \in \cM$ such that $uBu^* \subset C' \cap \cM$.
\end{definition}

One can check that if a separable II$_1$ factor $M$ has property $V$, then $M^{\omega}$ has property $\wV$, for any free ultrafilter $\omega$ on $\mathbb N$. 
Let $\Gamma$ be a countable group. Then $L(T_1(\Gamma))$ has property $V$ by \cite[Lemma 1]{MD69b}, hence  $L(T_1(\Gamma))^{\omega}$ has property $\wV$. On the other hand,  if $\Gamma$ is non-amenable, then we show that $L(T_0(\Gamma))^\omega$ does not have property $\wV$. More generally, we have the following theorem.

As in the previous section, we fix a countably cofinal ultrafilter $\cU$ on a set $S$. 

\begin{theorem}\label{uncountable}
Let $\Gamma$ be a non-amenable countable group. Let $k \geq 1$ and $\pmb{\alpha} \in \{0,1\}^k$. Let $M_{\pmb{\alpha}}(\Gamma)$ be as defined in Section \ref{sectionconstruction}.

Then $M_{\pmb{\alpha}}(\Gamma)^\cU$ has property $\wV$ at depth $k - 1$ (see Definition \ref{atdepthk}) if and only if $\alpha_{k} = 1$.

Moreover, if instead $\pmb\alpha$ is a sequence of $0$'s and $1$'s with length greater than $k$, then the same conclusion holds for arbitrary $\Gamma$ (possibly amenable).
\end{theorem}

Note that the moreover part of Theorem \ref{uncountable} follows from the first part of the statement. Indeed, if $k\geq 1$ and $\pmb\alpha$ has length greater that $k$, then Definition \ref{defintro} implies that the factor $M_{\pmb\alpha}(\Gamma)$ is of the form $M_{\pmb\beta}(\Lambda)$, for some non-amenable group $\Lambda$ and the truncated sequence $\pmb\beta = (\alpha_n)_{n=1}^k \in \{0,1\}^k$. 

Let us explain how this theorem implies our main results.

{\bf Proof of Theorem \ref{main}.} Assume that $M_{\pmb\alpha}(\Gamma)^{\cU}\simeq M_{\pmb\beta}(\Gamma)^{\cV}$, for two ultrafilters $\cU$ and $\cV$. If $\cU$ and $\cV$ are countably cofinal, then combining Theorem \ref{countable} and Theorem \ref{uncountable} leads to a contradiction. If one of the ultrafilters, say $\cU$, is not countably cofinal, then Lemma \ref{complete} readily implies that 
$M_{\pmb\alpha}(\Gamma)\simeq M_{\pmb\beta}(\Gamma)$, and \cite{MD69b} gives a contradiction. Alternatively, we may choose a free ultrafilter $\omega$ on $\mathbb N$ and derive a contradiction from $M_{\pmb\alpha}(\Gamma)^{\omega}\simeq M_{\pmb\beta}(\Gamma)^{\omega}$, as above.
\hfill$\blacksquare$

{\bf Proof of Corollary \ref{maincor}.}
Denote by $\tau_1$ the canonical trace on $C^*_r(T_{\pmb\alpha}(\Gamma))$, by $\tau_2$ the unique trace on $\cZ$, and let $\tau = \tau_1 \otimes \tau_2$. Note that the von Neumann algebra generated by $A_{\pmb\alpha}(\Gamma)$ acting via the GNS representation with respect to $\tau$ is $M_{\pmb\alpha}(\Gamma) \bar{\otimes} R \simeq M_{\pmb\alpha}(\Gamma)$. Indeed, $M_{\pmb\alpha}(\Gamma)$ is McDuff, as soon as $\pmb\alpha$ is non-empty.

It is an easy exercise  to show that for any ultrafilter $\cU$,  the von Neumann algebra generated by $A_{\pmb\alpha}(\Gamma)^\cU$ acting via the GNS representation associated with $\tau^\cU$ is precisely $M_{\pmb\alpha}(\Gamma)^\cU$.

Hence, in order to deduce Corollary \ref{maincor} from Theorem \ref{main}, we only need to check that $\tau^\cU$ is the only trace on $A_{\pmb\alpha}(\Gamma)^\cU$. 
Note that for any group $\Gamma$ and $\alpha \in \{0,1\}$, $T_\alpha(\Gamma)$ is equal to the increasing union of the free product groups $G_n:=\Sigma_n \ast \Lambda_{n+1}$, where $\Sigma_n$ is the subgroup generated by $\Lambda_1 \ast \dots \ast  \Lambda_n$ and $\Gamma_1 \oplus \dots \oplus \Gamma_n$ if $\alpha=0$ (respectively,  $(\Gamma_1 \oplus \dots \oplus \Gamma_n) \rtimes S_n$ if $\alpha = 1$). In particular, $T_{\pmb\alpha}(\Gamma)$ is an increasing union of Powers groups (see e.g.\ \cite{dlHS86}):
\begin{equation}\label{powers}
T_{\pmb\alpha}(\Gamma) = \bigcup_n G_n.
\end{equation}
Since $\cZ$ has a unique trace, \cite[Corollary 7]{dlHS86} implies that $A_{\pmb\alpha}(\Gamma)$ has a unique trace. Moreover, if $\Gamma$ is exact, so is $A_{\pmb\alpha}(\Gamma)$ and \cite[Theorem 8]{Oz13} shows that $A_{\pmb\alpha}(\Gamma)^\cU$ has the unique trace property. Let us treat the general case, when $\Gamma$ is not necessarily exact.

By \eqref{powers}, it is sufficient to show that for all families of integers $(k_s)_s$, the algebra $\prod_{\cU} (C^*_r(G_{k_s}) \otimes \cZ)$ has a unique trace. Since $\cZ^\cU$ has a unique trace and all the $G_{k_s}$'s are Powers groups, this is an easy consequence of \cite[Lemma 5]{dlHS86}.\hfill$\blacksquare$

The rest of this section is devoted to the prove Theorem \ref{uncountable}. As explained above we only need to prove the first statement. We will proceed by induction on $k$, and treat the base case and the inductive step in two separate subsections.

\subsection{The base case} 

The case $k = 0$ of Theorem \ref{uncountable} is dealt with by the following two lemmas.

\begin{lemma}\label{V1}
Let $\cU$ be a countably cofinal ultrafilter on a set $S$. Let $\G$ be a countable group and denote by $M = L(T_1(\Gamma))$. For every $s\in S$, let $t_s\geq 1$ be an integer. Let $\cM=\prod_{\cU}M^{{\bar\otimes} t_s}$. 

Then $\cM$ has property $\wV$.
\end{lemma}
{\it Proof.} Recall from Section \ref{sectionconstruction} that $T_1(\Gamma)$ is generated by $\widetilde\Gamma\rtimes S_{\infty}$ and $\Lambda_j, j\geq 1$, where  $\Gamma_i, i\geq 1$, and $\Lambda_j, j\geq 1,$ are isomorphic copies of $\Gamma$ and $\mathbb Z$, respectively, and  $S_{\infty}$ acts on $\widetilde\Gamma=\oplus_{i\geq1}\Gamma_i$  by permutations, with the only relations  that $\Gamma_i$ and $\Lambda_j$ commute whenever $i\geq j$.  Put $P=L(\widetilde\Gamma)$. 
Then Corollary \ref{combo} provides a separable subalgebra $A\subset\cM$ such that $A'\cap\cM\subset\prod_{\cU}P^{{\bar\otimes} t_s}$.

For $n' > n\geq 1$, recall that $\widetilde\Gamma_{n,n'}=\bigoplus_{n'>i\geq n}\Gamma_i$. For $n\geq 1$, let $H_n<T_1(\Gamma)$ be the subgroup generated by $\widetilde\Gamma_{1,n+1}\rtimes S_{n}$ and $\Lambda_1,...,\Lambda_n$, where we view $S_{n}$ as the group of all permutations of $\{1,2,...\}$ leaving each $k>n$ fixed.  Denote by  $R_n=L(\widetilde\Gamma_{1,n+1})$ and $M_n=L(H_n)$.

Let $B\subset\prod_{\cU}P^{{\bar\otimes} t_s}$ and $C\subset\cM$ be separable subalgebras. Since $\cup_{m\geq 1}R_m$ is weakly dense in $P$ and $\cup_{n\geq 1}M_n$ is weakly dense in $M$, then by Lemma \ref{ccofinal} exist integers 
$m_s,n_s\geq 1$, for $s\in S$, such that \begin{equation}\label{BC}B\subset\prod_{\cU}R_{m_s}^{{\bar\otimes} t_s}\;\;\text{and}\;\; C\subset\prod_{\cU}M_{n_s}^{{\bar\otimes} t_s}.\end{equation}

Finally, for every $s\in S$, let $\sigma_s\in S_{\infty}$ be a permutation such that $\sigma_s(k)>n_s$, for any $1\leq k\leq m_s$.
Then $\sigma_s\widetilde\Gamma_{1,m_s+1}\sigma_s^{-1}\subset\oplus_{i>n_s}\Gamma_i$, and hence $\sigma_s\widetilde\Gamma_{1,m_s+1}\sigma_s^{-1}$ commutes with $H_{n_s}$.  Thus, the unitary element $u_s=u_{\sigma_s}\in \ms U(M)$ satisfies $u_sR_{m_s}u_s^*\subset M_{n_s}'\cap M$. Therefore,  letting $u=(u_s^{\otimes t_s})_s\in\ms U(\cM)$, the equation \eqref{BC} implies that $uBu^*\subset C'\cap\cM$, as desired.
\hfill$\blacksquare$

\begin{lemma}\label{V2}
Let $\cU$ be an ultrafilter on a set $S$. Let $\G$ be a countable non-amenable group and denote by $M = L(T_0(\Gamma))$.
For every $s\in S$, let $Q_s$ be a tracial von Neumann algebra. 

Then any intermediate subalgebra $M^\cU \subset \cM \subset \prod_\cU (M {\bar\otimes} Q_s)$ does not have property $\wV$.
\end{lemma}

Lemma \ref{V2} strengthens \cite[Lemma 3]{MD69b}. In order to prove it, we will need two additional results.

Recall from Section \ref{sectionconstruction} that $T_0(\Gamma)$ is generated by $\widetilde\Gamma=\oplus_{i\geq 1}\Gamma_i$ and $\widetilde\Lambda=*_{j\geq 1}\Lambda_j$, where $\Gamma_i,\Lambda_j$ are isomorphic copies of $\Gamma,\mathbb Z$, respectively, with the only relations being that $\Gamma_i$ and $\Lambda_j$ commute whenever $i\geq j\geq 1$.
For $n\geq 1$, we denote by $\pi_n:\Gamma\rightarrow\widetilde\Gamma$ the canonical embedding with $\pi_n(\Gamma)=\Gamma_n$.
We also let $\widetilde\Gamma_n=\oplus_{i\geq n}\Gamma_i$.

\begin{lemma}\label{normal}
Let $g\in T_0(\Gamma), g'\in\widetilde\Gamma_{n+1}$, and  $g''\in\Gamma_{n}$, for some $n\geq 1$. Assume that $g'gg''=g$.

Then $g'=g''=e$. 
\end{lemma}

{\it Proof.}
Fix $n$, $g,g',g''$ as in the statement of the lemma satisfying $g'gg'' = g$.

Note that $T_0(\Gamma)$ splits as an amalgamated free product $T_0(\Gamma) = \Sigma \ast_{\widetilde \Gamma_{n+1}} \Delta$, where 
\begin{itemize}
\item $\Sigma < T_0(\Gamma)$ is the subgroup generated by $\widetilde \Gamma$, $\Lambda_1,..., \Lambda_n$;
\item $\Delta < T_0(\Gamma)$ is the subgroup generated by $\widetilde \Gamma_{n+1}$ and $\Lambda_{n+1}, \Lambda_{n+2}, ...$.
\end{itemize}
Then $g'' \in \Sigma$ is conjugate inside $T_0(\Gamma)$ to $g' \in \widetilde\Gamma_{n+1}$. By the General fact below, this implies that $g''$ and is actually conjugate inside $\Sigma$ to an element of $\widetilde\Gamma_{n+1}$. Now note that $\widetilde\Gamma_{n+1}$ is normal inside $\Sigma$ (it is even in product position). This forces $g''$ to belong to $\Gamma_n \cap \widetilde\Gamma_{n+1}$. Hence $g'' = e$, and further $g' = e$.

{\bf General fact.} Consider an amalgamated free product of groups $A = A_1 \ast_{A_0} A_2$. Assume that two elements $a_1 \in A_1$ and $a_2 \in A_2$ are conjugate inside $A$. Then $a_1$ is conjugate inside $A_1$ to an element in $A_0$.

Indeed, assume that $a_2 = ha_1h^{-1}$ for some $h \in A$. If $h \in A_1$ then clearly $a_2 \in A_1 \cap A_2$ and we are done. Otherwise, write $h$ as a product $h = h_0b$ for some $b \in A_1$ and some reduced word $h_0 \in A$ with rightmost letter in $A_2 \setminus A_0$. Then $ba_1b^{-1}$ lies inside $A_0$, for if not $h_0(ba_1b^{-1})h_0^{-1}$ would be a reduced word with a letter $bab^{-1} \in A_1 \setminus A_0$, so it could not lie inside $A_2$.
\hfill$\blacksquare$

\begin{lemma}\label{nonamen}
There exist $g_1,...,g_m\in\Gamma$ and $C>0$ such that the following holds:

For any $n\geq 1$, any unitaries $v_1,...,v_m\in\ms U(L(\widetilde\Gamma_{n+1}))$, and any $\xi\in M$ we have that 

$$\|\xi\|_2\leq C\sum_{k=1}^m\|u_{\pi_n(g_k)}\xi-\xi v_{k}\|_2.$$
\end{lemma}

{\it Proof.} Let $\lambda_n:\Gamma_n\rightarrow\mathcal U(\ell^2(T_0(\Gamma)/\widetilde\Gamma_{n+1}))$ be the quasi-regular representation of $\Gamma_n$.
Lemma \ref{normal} implies that $\Gamma_n$ acts freely on $T_0(\Gamma)/\widetilde\Gamma_{n+1}$. Thus,  $\lambda_n$ is a multiple of the left regular representation of $\Gamma_n$, hence $\lambda_n\circ\pi_n$ is a multiple of the left regular representation of $\Gamma$. 
 Since $\Gamma$ is non-amenable, there exist $g_1,...,g_m\in\Gamma$ and $C>0$ such that  for every $\xi\in\ell^2(T_0(\Gamma)/\widetilde\Gamma_{n+1})$ and $n\geq 1$ we have \begin{equation}\label{nonamen1}\|\xi\|_2\leq C\sum_{k=1}^m\|\lambda_n(\pi_n(g_k))\xi-\xi\|_2.\end{equation}
 
 We identify $L^2(M)\equiv\ell^2(T_0(\Gamma))$ as usual, via the unitary given by $u_g\mapsto\delta_g$, for any $g\in T_0(\Gamma)$. 
Fix $n\geq 1$. 
For $S\subset T_0(\Gamma)$, we denote by $P_S$ the orthogonal projection from $\ell^2(T_0(\Gamma))$ onto the  $\|\cdot\|_2$-closed linear span of $\{\delta_g\,|\,g\in S\}$. 
We define $T:\ell^2(T_0(\Gamma))\rightarrow \ell^2(T_0(\Gamma)/\widetilde\Gamma_{n+1})$ by letting $$T(\xi)(g\widetilde\Gamma_{n+1})=\|P_{g\widetilde\Gamma_{n+1}}(\xi)\|_2,\;\;\text{for every $\xi\in\ell^2(T_0(\Gamma))$ and $g\in T_0(\Gamma)$}.$$

Then for every $\xi,\eta\in\ell^2(T_0(\Gamma))$, $g\in\Gamma_n$, and $v \in \ms U(L(\widetilde\Gamma_{n+1}))$ we have that  \begin{equation}\label{nonamen2}
\|T(\xi)-T(\eta)\|_2\leq \|\xi-\eta\|_2,\;\;\;\;\;T(u_g\xi)=\lambda_n(g)(T(\xi)),\;\;\;\;\text{and}\;\;\;\;T(\xi u)=T(\xi).
\end{equation}
For the last identity, just notice that since the $\|\cdot\|_2$-closed linear span of $\{\delta_h\,|\,h\in g\widetilde\Gamma_{n+1}\}$ is a right $L(\widetilde\Gamma_{n+1})$-module, we have that  $P_{g\widetilde\Gamma_{n+1}}(\xi v)=P_{g\widetilde\Gamma_{n+1}}(\xi)v$, for every $g\in T_0(\Gamma)$ and $v \in \ms U(L(\widetilde\Gamma_{n+1}))$.  

Finally, let $\xi\in\ell^2(T_0(\Gamma))$ and $v_1,...,v_m\in\ms U(L(\widetilde\Gamma_{n+1}))$. By applying \eqref{nonamen1} to $T(\xi)$ and using \eqref{nonamen2} we deduce
\[\|\xi\|_2=\|T(\xi)\|_2 \leq C\sum_{k=1}^m\|T(u_{\pi_n(g_k)}\xi)-T(\xi v_k)\|_2 \leq C\sum_{k=1}^m\|u_{\pi_n(g_k)}\xi-\xi v_k\|_2.\]
This finishes the proof.
 \hfill$\blacksquare$

\subsection*{Proof of Lemma \ref{V2}} Assume by contradiction that there exists an intermediate subalgebra $M^{\cU}\subset\cM\subset\prod_{\cU}(M{\bar\otimes} Q_s)$ with property $\wV$. Let $A\subset\cM$ be as in Definition \ref{propV}. For $n\geq 1$, put $P_n=L(\widetilde\Gamma_n)$.
Applying Corollary \ref{combo2}, we can find integers $i_s\geq 1$, for every $s\in S$, such that 
\[\prod_{\cU}P_{i_s}\subset A'\cap\cM.\]
Let $\rho:\Gamma\rightarrow\mathcal U(\prod_{\cU}P_{i_s})$  be the homomorphism given by $\rho(g)=(u_{\pi_{i_s}(g)})_s$, for $g\in\Gamma$.
Consider $B\subset A' \cap \cM$ the (separable) subalgebra generated by $\rho(\Gamma)$. Also, use Corollary \ref{combo} to get a separable subalgebra $C\subset M^{\cU}$ such that 
\[C'\cap\cM\subset\prod_{\cU}(P_{i_s+1}{\bar\otimes} Q_s).\]
Since $\cM$ has property $\wV$, there exists a unitary $v \in \cM \subset \prod_{\cU}(M{\bar\otimes} Q_s)$ such that 
\[v^*Bv\subset C' \cap \cM \subset  \prod_{\cU}(P_{i_s+1}{\bar\otimes} Q_s).\]
 
 Represent $v=(v_s)_s$, where $v_s\in M{\bar\otimes} Q_s$ is a unitary, for any $s\in S$. Let $g_1,...,g_m\in\Gamma$ and $C>0$ be given by Lemma 
 \ref{nonamen}. For $1\leq k\leq m$, denote by $u_k=v^*\rho(g_k)v$. Since $u_k\in\prod_{\cU}(P_{i_s+1}{\bar\otimes} Q_s)$, we can represent $u_k=(u_{k,s})_s$, where $u_{k,s}\in P_{i_s+1}{\bar\otimes} Q_s$ is a unitary. 
 Since $\rho(g_k)v=vu_k$, we get
 
 $$\lim\limits_{\cU}\|u_{\pi_{i_s}(g_k)}v_s-v_su_{k,s}\|_2=0,\;\;\text{for every $1\leq k\leq m$}.$$
 
 Since $\|v_s\|_2=1$, for every $s\in S$, this clearly contradicts the conclusion of Lemma \ref{nonamen}. \hfill$\blacksquare$

\subsection{The inductive step} 

Theorem \ref{uncountable} clearly follows from the next two lemmas.

\begin{lemma}
Let $\Gamma$ be a countable group. Use the notations of Section \ref{sectionconstruction}. 

For any $k \geq 1$ and $\pmb{\alpha} \in \{0,1\}^k$ such that $\alpha_k = 1$, any integers $t_s\geq 1, s\in S$, we have that $\prod_\cU {M_{\pmb\alpha}(\Gamma)}^{{\bar{\otimes}} t_s}$ has property $\wV$ at depth $k - 1$.
\end{lemma}

{\it Proof.} 
We proceed by induction on $k$. If $k = 1$, then $M_{\pmb\alpha}(\Gamma) =L(T_1(\Gamma))$ and the conclusion follows from Lemma \ref{V1}.

Assume the conclusion holds for some $k \geq 1$. Let $\pmb\alpha \in \{0,1\}^{k+1}$ be such that $\alpha_{k+1} = 1$. Then the sequence $\pmb\beta = (\alpha_{n+1})_{n = 1}^k$ has length $k$ and $\beta_{k} = \alpha_{k+1}=1$. Moreover,  we have that $M_{\pmb\alpha}(\Gamma) =L(T_{\alpha_1} (T_{\pmb\beta}(\Gamma)))$. Let $t_s\geq 1, s\in S$, be integers and denote by $\cM = \prod_{\cU} {M_{\pmb\alpha}(\Gamma)}^{{\bar{\otimes}} t_s}$.

By applying Lemma \ref{MDex} (to $T_{\pmb\beta}(\Gamma)$ and $\alpha_1$ instead of $\Gamma$ and $\alpha$)  we get that
 $A_i = \prod_\cU {P_{i_s}}^{{\bar{\otimes}} t_s}$, where $i = (i_s)_s\in I$, is a residual net for $\cM$.
 Since $\beta_{k} = 1$ and  $P_{m,m'}\cong {M_{\pmb\beta}(\Gamma)}^{{\bar{\otimes}} (m'-m)}$, for any $m'>m\geq 1$, then using the inductive hypothesis and repeating the end of the proof of Lemma \ref{count1} it follows that, at depth $1$, $\cM$ has property $\wV$ at depth $k - 1$. 
 
 This shows that $\cM$ has property $\wV$ at depth $k$, and finishes the proof.
\hfill$\blacksquare$

\begin{lemma}
Let $\Gamma$ be a non-amenable countable group and keep the notations from Section \ref{sectionconstruction}. 

For any $k \geq 1$ and $\pmb{\alpha} \in \{0,1\}^k$ such that $\alpha_k = 0$ and any family of tracial von Neumann algebras $Q_s\geq 1, s\in S$, no intermediate von Neumann subalgebra ${M_{\pmb\alpha}(\Gamma)}^\cU = \prod_\cU M_{\pmb\alpha}(\Gamma) \subseteq \cM \subseteq \prod_\cU (M_{\pmb\alpha}(\Gamma) {\bar\otimes} Q_s)$ has property $\wV$ at depth $k - 1$.
\end{lemma}
{\it Proof.} 
We proceed by induction on $k$. If $k = 1$, then $M_{\pmb\alpha}(\Gamma) = L(T_0(\Gamma))$ and since $\Gamma$ is non-amenable, the conclusion follows from Lemma \ref{V2}.

Now assume the conclusion holds for some $k \geq 1$. Let $\pmb\alpha \in \{0,1\}^{k+1}$ be such that $\alpha_{k+1} = 0$.  Then the sequence $\pmb\beta = (\alpha_{n+1})_{n = 1}^k$ has length $k$ and $\beta_{k} = \alpha_{k+1} = 0$. Moreover, $M_{\pmb\alpha}(\Gamma) =L(T_{\alpha_1} (T_{\pmb\beta}(\Gamma)))$.

Suppose by contradiction that there exist tracial von Neumann algebras $Q_s, s\in S$, and an algebra $\cM$ satisfying 
\[{M_{\pmb{\alpha}}(\Gamma)}^\cU \subset \cM \subset \prod_\cU (M_{\pmb{\alpha}}(\Gamma) {\bar\otimes} Q_s),\]
that has property $\wV$ at depth $k$.

Then, at depth $1$, $\cM$ has property $\wV$ at depth $k-1$.  Let $(A_i \subset B_i)_{i \in I}$ be the residual pair for the inclusion $M_{\pmb\alpha}(\Gamma)^\cU \subset \prod_\cU (M_{\pmb\alpha}(\Gamma) {\bar\otimes} Q_s)$ described in Lemma \ref{MDex2}, with $T_{\pmb\beta}(\Gamma)$ and $\alpha_1$ instead of $\Gamma$ and $\alpha$.  Since $\beta_{k+1} = 0$ and $P_{m,m'}\cong M_{\pmb\beta}(\Gamma)^{{\bar\otimes} (m'-m)}$, for any $m'>m\geq 1$, by using the inductive hypothesis and repeating the end of the proof of Lemma \ref{count2} we get a contradiction.
\hfill$\blacksquare$

\bibliographystyle{alpha1}

\end{document}